\documentclass[12pt,reqno]{amsart}

\usepackage{mathdots}

\usepackage{color}

\usepackage[dvipdfmx]{graphicx}

\usepackage{amssymb}

\usepackage{mathrsfs}

\DeclareMathAlphabet{\mathpzc}{OT1}{pzc}{m}{it}

\usepackage[all]{xy}

\usepackage{tikz}
\usetikzlibrary{arrows,decorations.pathmorphing,decorations.pathreplacing,positioning,shapes.geometric,shapes.misc,decorations.markings,decorations.fractals,calc,patterns}

\usepackage{float}

\usepackage[bottom]{footmisc}

\entrymodifiers={+!!<0pt,\fontdimen22\textfont2>}

\def\sA{\mathsf{A}}
\def\sB{\mathsf{B}}

\def\sK{\mathsf{K}}

\def\sT{\mathsf{T}}
\def\sU{\mathsf{U}}
\def\sV{\mathsf{V}}

\def\sX{\mathsf{X}}
\def\sY{\mathsf{Y}}
\def\sZ{\mathsf{Z}}

\def\ac{\operatorname{ac}}

\def\adots{\mathinner{\mkern1mu\raise1.0pt\vbox{\kern7.0pt\hbox{.}}\mkern2mu\raise4.0pt\hbox{.}\mkern2mu\raise7.0pt\hbox{.}\mkern1mu}}

\def\b{\operatorname{b}}

\def\cone{\operatorname{cone}}

\def\dddots{\mathinner{\mkern1mu\raise10.0pt\vbox{\kern7.0pt\hbox{.}}\mkern2mu\raise5.3pt\hbox{.}\mkern2mu\raise1.0pt\hbox{.}\mkern1mu}}
\def\dddotssmall{\mathinner{\mkern1mu\raise7.0pt\vbox{\kern7.0pt\hbox{.}}\mkern-1mu\raise4pt\hbox{.}\mkern-1mu\raise1.0pt\hbox{.}\mkern1mu}}

\def\Hom{\operatorname{Hom}}

\def\Prj{\operatorname{Prj}}

\def\s{\operatorname{s}}

\def\SL2{\operatorname{SL}_2}

\def\tac{\operatorname{tac}}

\numberwithin{equation}{section}

\renewcommand{\labelenumi}{(\roman{enumi})}

\newtheorem{Lemma}{Lemma}[section]
\newtheorem{Theorem}[Lemma]{Theorem}

\newtheorem{Corollary}[Lemma]{Corollary}

\theoremstyle{definition}

\newtheorem{Remark}[Lemma]{Remark}

\newtheorem*{bfhpg*}{}

\begin{document}

\setlength{\parindent}{0pt}
\setlength{\parskip}{7pt}

\title[Extensions, t-structures, triangles of
recollements]{Triangulated subcategories of extensions, stable
  t-structures, and triangles of recollements}

\author{Peter J\o rgensen}
\address{School of Mathematics and Statistics,
Newcastle University, Newcastle upon Tyne NE1 7RU, United Kingdom}
\email{peter.jorgensen@ncl.ac.uk}
\urladdr{http://www.staff.ncl.ac.uk/peter.jorgensen}

\author{Kiriko Kato}
\address{Department of Mathematics and Information Sciences, Osaka
Prefecture University, Osaka, Japan}
\email{kiriko@mi.s.osakafu-u.ac.jp}
\urladdr{http://www.mi.s.osakafu-u.ac.jp/\~{ }kiriko}




\begin{abstract} 

In a triangulated category $\sT$ with a pair of triangulated
subcategories $\sX$ and $\sY$, one may consider the subcategory of
extensions $\sX * \sY$.

We give conditions for $\sX * \sY$ to be triangulated and use them
to provide tools for constructing stable t-structures.  In
 particular, we show how to construct so-called triangles of
recollements, that is, triples of stable t-structures of the form $(
\sX,\sY )$, $( \sY,\sZ )$, $( \sZ,\sX )$.  We easily recover some
 triangles of recollements known from the literature.

\end{abstract}

\maketitle

\setcounter{section}{-1}
\section{Introduction}
\label{sec:introduction}

Let $\sT$ be a triangulated category and let $\sX,\sY \subseteq \sT$
be subcategories.  The subcategory of extensions is defined by
\[
  \sX * \sY
  = \bigg\{\, e \in \sT \,\bigg|
      \begin{array}{ll}
         \mbox{ there is a distinguished triangle } \\[2mm]
         \mbox{ $x \rightarrow e \rightarrow y$ in $\sT$
                with $x \in \sX$, $y \in \sY$}
      \end{array}
    \,\bigg\}.
\]
This is a classic object which, as far as we know, was introduced by
Beilinson, Bernstein, and Deligne in \cite[1.3.9]{BBD} and has been
used extensively since then, see for instance Bondal--Van den Bergh
\cite[sec.\ 2]{BVdB}, Iyama--Yoshino \cite[sec.\ 2]{IY}, and Rouquier
\cite[sec.\ 3]{R}. 

It has been well-known that $\sX * \sY$ is triangulated if 
$\Hom_{ \sT }( \sX,\sY ) = 0$. 
We have examined this condition to see that 
a substantial generalization provides exact characterization of 
$\sX * \sY$ to be triangulated. 


\bigskip
{\bf Theorem A. }
{\em
Let $\sX, \sY \subseteq \sT$ be triangulated subcategories.  
The subcategory $\sX * \sY$ is triangulated if and only if 
$\Hom_{ \sT / \sX \cap \sY }( Q\sX , Q\sY ) = 0$, where
  $\sT / \sX \cap \sY$ is the Verdier quotient and $Q : \sT
  \rightarrow \sT / \sX \cap \sY$ is the quotient functor.

}
\bigskip

Note that throughout, all subcategories
are full and closed under isomorphisms in the ambient category, and
that if $Q : \sT \rightarrow \sT'$ is a functor and $\sX \subseteq
\sT$ a subcategory, then $Q\sX$ denotes the isomorphism closure of
$\{\, Qx \,|\, x \in \sX \,\}$ in $\sT'$, viewed as a subcategory.

Theorem A is a main ingredient in the proof of Theorem B below which
is a tool for constructing stable t-structures.  Recall that a stable
t-structure in $\sT$ is a pair $( \sX,\sY )$ of subcategories which
are stable under $\Sigma$ and satisfy $\Hom_{ \sT }( \sX,\sY ) = 0$
and $\sX * \sY = \sT$, see \cite[def.\ 9.14]{M}.  Stable t-structures
are important in several settings which involve triangulated
categories, see for instance \cite{IKM}.  In particular, they play a
prominent role in algebraic geometry where they are known as
semi-orthogonal decompositions, see for instance \cite{H}.

Note that if $\sX, \ldots, \sY \subseteq \sT$ are subcategories, then
$\langle \sX, \ldots, \sY \rangle$ denotes the smallest triangulated
subcategory containing them.

\bigskip
{\bf Theorem B. }
{\em
Each row in the following table defines a triangulated subcategory
$\sU \subseteq \sT$ and considers the quotient functor $Q : \sT
\rightarrow \sT / \sU$.  It then shows one or more stable t-structures
in $\sT / \sU$.
}
\[
\mbox{
\bgroup
\def\arraystretch{2.5}
\begin{tabular}{c|c|c|c}
  \begin{minipage}[c]{3.75cm} Given triangulated \\ subcategories of $\sT$... \end{minipage} & Which satisfy... & Define... & \begin{minipage}[c]{4.0cm} Then there are stable \\ t-structures in $\sT/\sU$... \end{minipage} \\ \cline{1-4}
  $\sX$, $\sY$ & \;\;\;\;\; \begin{minipage}[c]{2.5cm}$\sX * \sY = \sT$ \end{minipage}& $\sU = \sX \cap \sY$ & \begin{minipage}[c]{2.5cm}$( Q\sX , Q\sY )$ \end{minipage}\\[2mm] \cline{1-4}
  $\sX$, $\sY$, $\sZ$ & \;\;\;\;\;\;\begin{minipage}[c]{2.5cm} \vspace{1ex}$\sX * \sY = \sT$ \\[2mm] $\sY * \sZ = \sT$ \end{minipage}& $\sU = \langle \sX \cap \sY , \sY \cap \sZ \rangle$ & \begin{minipage}[c]{2.5cm} \vspace{1ex}$( Q\sX , Q\sY )$ \\[2mm] $( Q\sY , Q\sZ )$ \end{minipage} \\[4mm] \cline{1-4}
  $\sX$, $\sY$, $\sZ$ & \;\;\;\;\;\:\begin{minipage}[c]{2.5cm} \vspace{1ex}$\sX * \sY = \sT$ \\[2mm] $\sY * \sZ = \sT$ \\[2mm] $\sZ * \sX = \sT$ \end{minipage}& $\sU = \langle \sX \cap \sY , \sY \cap \sZ , \sZ \cap \sX \rangle$ & \begin{minipage}[c]{2.5cm} \vspace{1ex}$( Q\sX , Q\sY )$ \\[2mm] $( Q\sY , Q\sZ )$ \\[2mm] $( Q\sZ , Q\sX )$ \end{minipage} \\
\end{tabular}
\egroup
}
\]
\bigskip

Note that the penultimate row of the table gives two ``adjacent''
stable t-structures, and that this is equivalent to giving a
recollement, see \cite[prop.\ 1.8]{IKM}.  The last row gives three
``cyclically adjacent'' stable t-structures.  In the terminology of
\cite[def.\ 1.9]{IKM}, the triple $( Q\sX, Q\sY, Q\sZ )$ is a
``triangle of recollements''.  It is somewhat surprising that such a
structure is possible, but it provides a setup with a large number of
pleasant symmetries, see for instance \cite[props.\ 1.10 and
1.16]{IKM}.  We show in Section \ref{sec:examples} how Theorem B
easily recovers two triangles of recollements known from \cite[thm.\
5.8]{IKM} and \cite[thm.\ 2.10]{JK}.

Our third main theorem is a variation on a part of Theorem B; it is
also proved by means of Theorem A.  Note that a thick subcategory is a
triangulated subcategory closed under direct summands.  If $( \sX,\sY
)$ is a stable t-structure, then $\sX$ and $\sY$ are thick by
\cite[prop.\ 9.15]{M}.

\bigskip
{\bf Theorem C. }
{\em
Let $\sU \subseteq \sT$ be a thick subcategory, $Q : \sT \rightarrow
\sT / \sU$ the quotient functor.  There is a bijection 
}
\[
  \bigg\{\, ( \sX , \sY )
  \,\bigg|
    \begin{array}{l}
      \mbox{ $\sX, \sY \subseteq \sT$ thick subcategories } \\[2mm]
      \mbox{ with $\sX \cap \sY = \sU$ and $\sX * \sY = \sT$ }
    \end{array}
  \,\bigg\}
  \leftrightarrow
  \{\, \mbox{ stable t-structures in $\sT/\sU$ } \,\}
\]
{\em given by $( \sX,\sY ) \mapsto ( Q\sX,Q\sY )$. }
\bigskip

The paper is organised as follows: Section \ref{sec:proof_of_A} gives
the proof of Theorem A.  Section \ref{sec:lemmas} proves some lemmas
needed later.  Sections \ref{sec:proof_of_B} and \ref{sec:proof_of_C}
give the proofs of Theorems B and C.  Section \ref{sec:another} gives
another result based on Theorem A, and Section \ref{sec:examples} uses
Theorem B to recover two triangles of recollements known from the
literature.

\section{Proof of Theorem A}
\label{sec:proof_of_A}

We shall begin with the following observation 
which helps us to prove Theorem A. 

\begin{Lemma}\label{lem:factorization}
Let $\sX, \sY \subseteq \sT$ be triangulated subcategories.  The
following conditions are e\-qui\-va\-lent.
\begin{enumerate}

  \item  $\sX * \sY$ is triangulated.

\smallskip

  \item  $\sY * \sX \subseteq \sX * \sY$.

\smallskip

  \item  Each morphism $x \rightarrow y$ with $x \in \sX$, $y \in \sY$
    factors through an object from $\sX \cap \sY$.

\end{enumerate}

\end{Lemma}
\bigskip

(i) $\Leftrightarrow$ (ii): (i) $\Rightarrow$ (ii) is obvious. 
(ii)$\Rightarrow$(i): Since $\sY * \sX \subseteq \sX * \sY$, we can
compute as follows using that $*$ is associative by \cite[lem.\
1.3.10]{BBD}.
\[
  ( \sX * \sY ) * ( \sX * \sY )
  = \sX * ( \sY * \sX ) * \sY
  \subseteq \sX * ( \sX * \sY ) * \sY
  = ( \sX * \sX ) * ( \sY * \sY )
  = \sX * \sY.
\]
So $\sX * \sY$ is closed under extensions, and since it is also closed
under $\Sigma^{ \pm 1 }$, it is triangulated. 

We should mention (i) $\Leftrightarrow$ (ii) immediately follows also from 
Lemma~\ref{lem:angles0} since
(i) equivalently says $X*Y = \langle X, Y \rangle$.

(ii)$\Rightarrow$(iii): Let $x \rightarrow y$ be a morphism with $x
\in \sX$, $y \in \sY$.  Completing to a distinguished triangle $x
\rightarrow y \rightarrow c$, we have $c \in \sY * \sX \subseteq \sX *
\sY$, so there is a distinguished triangle $x_1 \rightarrow c
\rightarrow y_1$ with $x_1 \in \sX$, $y_1 \in \sY$.  The octahedral
axiom provides a commutative diagram where rows and columns are
distinguished triangles.
\[
  \xymatrix {
    \Sigma^{-1}x_1 \ar[r] \ar[d] & 0 \ar[r] \ar[d] & x_1 \ar[d] \\
    x \ar[r] \ar[d] & y \ar[r] \ar@{=}[d] & c \ar[d] \\
    w \ar[r] & y \ar[r] & y_1 \\
            }
\]
The first column shows $w \in \sX$ and the last row shows $w \in
\sY$, so $x \rightarrow w \rightarrow y$ is a factorization as
claimed in (iii).

(iii)$\Rightarrow$(ii): If $e \in \sY * \sX$ is given then there is a
distinguished triangle $y \rightarrow e \rightarrow x \stackrel{ \xi
}{ \rightarrow } \Sigma y$ with $y \in \sY$, $x \in \sX$.  Use (iii)
to factorise $\xi$ as $x \rightarrow w \rightarrow \Sigma y$ with $w
\in \sX \cap \sY$.  The octahedral axiom again gives a commutative
diagram where rows and columns are distinguished triangles.
\[
  \xymatrix {
    x_1 \ar@{=}[r] \ar[d] & x_1 \ar[r] \ar[d] & 0 \ar[d] \\
    e \ar[r] \ar[d] & x \ar^{\xi}[r] \ar[d] & \Sigma y \ar@{=}[d] \\
    y_1 \ar[r] & w \ar[r] & \Sigma y \\
            }
\]
The second column shows $x_1 \in \sX$ and the last row shows $y_1 \in
\sY$, so the first column shows $e \in \sX * \sY$.

\hfill $\Box$

\begin{Remark} 
The equivalence between (i) and (ii) is 
a triangulated analogue of \cite[thm.\ 3.2]{Y} 
on Serre subcategories of module categories. 
\end{Remark}

{\em Proof} of Theorem A. It suffices to show that 
(iii) of Lemma~\ref{lem:factorization} is equivalent to 
the following. 

\begin{enumerate}
\setcounter{enumi}{3}
\item 
  We have $\Hom_{ \sT / \sX \cap \sY }( Q\sX , Q\sY ) = 0$, where
  $\sT / \sX \cap \sY$ is the Verdier quotient and $Q : \sT
  \rightarrow \sT / \sX \cap \sY$ is the quotient functor.
\end{enumerate}

Indeed, this is due to Verdier.  Namely, using
his terminology, condition (iii) is equivalent to $\sY$ being
``$\sX$-localisante \`a droite'' by \cite[prop.\ II.2.3.5(a)]{V}.  By
\cite[prop.\ II.2.3.5(c)]{V} this is equivalent to $\sY / \sX \cap
\sY$ consisting of objects which are ``$R$-libre  \`a droite'' where $R :
\sT / \sX \cap \sY \rightarrow \sT / \sX$ is the quotient functor.
Finally, this is equivalent to condition (iv) by \cite[prop.\
II.2.3.3(a)]{V}.

\hfill $\Box$

\begin{Remark}
We could have added the following conditions to Lemma~\ref{lem:factorization} 
\begin{enumerate}
\setcounter{enumi}{4}

  \item  If $\sX' \subseteq \sT$ is a triangulated subcategory such
    that $\sX \cap \sY \subseteq \sX' \subseteq \sX$, then $\sX' *
    \sY$ is triangulated.

\smallskip

  \item  If $\sY' \subseteq \sT$ is a triangulated subcategory such
    that $\sX \cap \sY \subseteq \sY' \subseteq \sY$, then $\sX *
    \sY'$ is triangulated.

\end{enumerate}
Namely, there are implications as follows between these and conditions
(i)--(iii) of Lemma~\ref{lem:factorization}  and (iv) in the proof of Theorem A.

(iii)$\Rightarrow$(v): 
It is enough to see $\sY * \sX' \subseteq \sX' * \sY$, but this
follows by a small variation of the proof of (iii)$\Rightarrow$(ii) in
Lemma~\ref{lem:factorization}.

(iii)$\Rightarrow$(vi) is similar, and (v)$\Rightarrow$(i) and
(vi)$\Rightarrow$(i) are trivial.
\end{Remark}

The following consequence of Theorem A is well-known, but we do not
know a published reference.

\begin{Corollary}
\label{cor:A}
If $\sX, \sY \subseteq \sT$ are triangulated subcategories with
$\Hom_{ \sT }( \sX,\sY ) = 0$, then $\sX * \sY$ is triangulated.
\end{Corollary}

\begin{proof}
Since $\Hom_{ \sT }( \sX,\sY ) = 0$ implies that  $\sX \cap \sY$ is trivial,  
$Q : \sT \rightarrow \sT / \sX \cap \sY$ is an equivalence. 
Hence  $\Hom_{ \sT / \sX \cap \sY }( Q\sX , Q\sY ) =\Hom_{ \sT }( \sX,\sY ) = 0$. 
Due to Theorem A, $\sX * \sY$ is triangulated.
\end{proof}

\section{Lemmas}
\label{sec:lemmas}

The following lemma uses the notation $\sX^{ *n } = \sX * \cdots *
\sX$ with $\sX$ repeated $n$ times.

\begin{Lemma}
\label{lem:angles0}
If $\sX, \sY \subseteq \sT$ are subcategories closed under $\Sigma^{
  \pm 1 }$ and containing $0$, then 
\[
  \langle \sX,\sY \rangle = \bigcup_{ n \geq 1 } ( \sX * \sY )^{*n}.
\]
\end{Lemma}

\begin{proof}
The inclusion $\supseteq$ is clear and $\subseteq$ holds because the
right hand side contains $\sX$ and $\sY$, and is closed under
$\Sigma^{ \pm 1 }$ and extensions so is triangulated.
\end{proof}

\begin{Lemma}
\label{lem:angles}
Let $\sX' \subseteq \sX \subseteq \sT$ and $\sY' \subseteq \sY
\subseteq \sT$ be triangulated subcategories such that $\sX * \sY$ is
triangulated.  Then
\begin{enumerate}

  \item  $\sX \cap \langle \sX',\sY \rangle = \langle \sX',\sX \cap \sY \rangle$,

\smallskip

  \item  $\langle \sX,\sY' \rangle \cap \sY = \langle \sX \cap \sY,\sY' \rangle$.

\end{enumerate}
\end{Lemma}

\begin{proof}
Part (ii) is clear from part (i) by symmetry.  In (i) the inclusion
$\supseteq$ is clear. 

To prove $\subseteq$ in (i), it is enough by Lemma
\ref{lem:angles0} to show $\sX \cap ( \sY * \sX' )^{*n} \subseteq
\langle \sX',\sX \cap \sY \rangle$ for each $n \geq 1$.  We do so by
showing the stronger statement
\[
  Q_n
  = \left\{
      \begin{array}{ll}
        \mbox{(a)} & \sX \cap ( \sY * \sX' )^{*n} \subseteq \langle \sX',\sX \cap \sY \rangle, \\[2mm]
        \mbox{(b)} & \sX \cap \big( ( \sY * \sX' )^{*n} * \sY \big) \subseteq \langle \sX',\sX \cap \sY \rangle \\
      \end{array}
    \right.
\]
for each $n \geq 1$ by induction.

$n=1$: To show $Q_1$, part (a), let $e \in \sX \cap ( \sY * \sX' )$ be
given.  There is a distinguished triangle $y \rightarrow e \rightarrow
x'$ with $y \in \sY$, $x' \in \sX'$.  Since $e$ and $x'$ are both in
$\sX$, so is $y$, so we have $y \in \sX \cap \sY$ whence the
distinguished triangle shows $e \in ( \sX \cap \sY ) * \sX' \subseteq
\langle \sX',\sX \cap \sY \rangle$.

To show $Q_1$, part (b), let $e \in \sX \cap \big( ( \sY * \sX' ) *
\sY \big)$ be given.  Since $*$ is associative we have $e \in \sY * (
\sX' * \sY )$ so there is a distinguished triangle $y \rightarrow e
\rightarrow a$ with $y \in \sY$, $a \in \sX' * \sY$.  There is also a
distinguished triangle $x'_1 \rightarrow a \rightarrow y_1$ with $x'_1
\in \sX'$, $y_1 \in \sY$.

We have $e \in \sX$, $y_1 \in \sY$, so by Lemma \ref{lem:factorization} 
the composition $e
\rightarrow a \rightarrow y_1$ factors as $e \rightarrow w \rightarrow
y_1$ for an object $w \in \sX \cap \sY$.  The $9$-lemma \cite[prop.\
4.9]{M} now provides a commutative diagram where rows and columns are
distinguished triangles.
\[
  \xymatrix {
    y_3 \ar[r] \ar[d] & x_2 \ar[r] \ar[d] & x'_1 \ar[d] \\
    y \ar[r] \ar[d] & e \ar[r] \ar[d] & a \ar[d] \\
    y_2 \ar[r] & w \ar[r] & y_1 \\
            }
\]
The second column shows $x_2 \in \sX$.  The third row shows $y_2 \in
\sY$, whence the first column shows $y_3 \in \sY$ and the first row
shows $x_2 \in \sY * \sX'$.  Hence $x_2 \in \sX \cap ( \sY * \sX' )$
whence $Q_1$, part (a) shows $x_2 \in \langle \sX',\sX \cap \sY
\rangle$.  Then the second column shows $e \in \langle \sX',\sX \cap
\sY \rangle$ as desired.

$n \geq 2$: To show $Q_n$, part (a), let $e \in \sX \cap ( \sY * \sX'
)^{*n}$ be given.  Since $*$ is associative there is a distinguished
triangle $y \rightarrow e \rightarrow a$ with $y \in \sY$, $a \in (
\sX' * \sY )^{*(n-1)} * \sX'$.  There is also a distinguished triangle
$b \rightarrow a \rightarrow x'$ with $b \in ( \sX' * \sY )^{*(n-1)}$,
$x' \in \sX'$.  The octahedral axiom provides a commutative diagram
where rows and columns are distinguished triangles.
\[
  \xymatrix {
    y \ar[r] \ar@{=}[d] & x \ar[r] \ar[d] & b \ar[d] \\
    y \ar[r] \ar[d] & e \ar[r] \ar[d] & a \ar[d] \\
    0 \ar[r] & x' \ar@{=}[r] & x' \\
            }
\]
The second column shows $x \in \sX$ and the first row shows $x \in \sY
* ( \sX' * \sY )^{ *(n-1) } = ( \sY * \sX' )^{ *(n-1) } * \sY$.
Statement $Q_{n-1}$, part (b), hence shows $x \in \langle \sX',\sX
\cap \sY \rangle$, and then the second column shows $e \in \langle
\sX',\sX \cap \sY \rangle$ as desired.

To show $Q_n$, part (b), let $e \in \sX \cap \big( ( \sY * \sX' )^{*n}
* \sY \big)$ be given.  Since $*$ is associative there is a
distinguished triangle $y \rightarrow e \rightarrow a$ with $y \in
\sY$, $a \in ( \sX' * \sY )^{ *n }$.  There is also a distinguished
triangle $b \rightarrow a \rightarrow y_1$ with $b \in \sX' * ( \sY *
\sX' )^{ *(n-1) }$, $y_1 \in \sY$.  

As above, the composition $e \rightarrow a \rightarrow y_1$ factors as
$e \rightarrow w \rightarrow y_1$ for an object $w \in \sX \cap \sY$,
and we obtain another commutative diagram where rows and columns are
distinguished triangles.
\[
  \xymatrix {
    y_3 \ar[r] \ar[d] & x \ar[r] \ar[d] & b \ar[d] \\
    y \ar[r] \ar[d] & e \ar[r] \ar[d] & a \ar[d] \\
    y_2 \ar[r] & w \ar[r] & y_1 \\
            }
\]
The second column shows $x \in \sX$.  The third row shows $y_2 \in
\sY$ whence the first column shows $y_3 \in \sY$ and the first row
shows $x \in ( \sY * \sX' )^{*n}$.  Hence $x \in \sX \cap ( \sY * \sX'
)^{*n}$ whence $Q_{ n }$, part (a) shows $x \in \langle \sX',\sX \cap
\sY \rangle$.  The second column shows $e \in \langle \sX',\sX \cap
\sY \rangle$ as desired.
\end{proof}

If $Q : \sT \rightarrow \sT'$ is a functor and $\sX' \subseteq \sT'$
is a subcategory then $Q^{-1}\sX' = \{\, t \in \sT \,|\, Qt \in \sX'
\,\}$ will be viewed as a subcategory of $\sX$.

\begin{Remark}
\label{rmk:Verdier}
If $\sU \subseteq \sT$ is a thick subcategory, then $\sX \mapsto
Q\sX$ gives a bijection
\[
  \bigg\{\,
    \begin{array}{l}
      \mbox{ triangulated subcategories $\sX \subseteq \sT$ } \\[2mm]
      \mbox{ with $\sU \subseteq \sX \subseteq \sT$ }
    \end{array}
  \,\bigg\}
  \leftrightarrow
  \{\, \mbox{triangulated subcategories of $\sT / \sU$} \,\}
\]
with inverse $\sY \mapsto Q^{ -1 }\sY$.  Under the bijection, thick
subcategories correspond to thick subcategories.  See \cite[prop.\
II.2.3.1]{V}.
\end{Remark}

\begin{Lemma}
\label{lem:Q}
Let $\sU \subseteq \sT$ be a thick subcategory and let $Q : \sT
\rightarrow \sT / \sU$ be the quotient functor.
\begin{enumerate}

  \item  If $\sX, \sY \subseteq \sT$ are triangulated subcategories
with $\sU \subseteq \sX \cap \sY$, then

\smallskip

\begin{enumerate}

  \item  $Q( \sX \cap \sY ) = Q\sX \cap Q\sY$,

\smallskip

  \item  $Q( \sX * \sY ) = Q\sX * Q\sY$. 

\end{enumerate}

\smallskip

  \item  If $\sX \subseteq \sT$ is any subcategory then
$Q^{-1}Q\sX = \sU * \sX * \sU$.

\end{enumerate}
\end{Lemma}

\begin{proof}
(i.a): The inclusion $\subseteq$ is clear.  To see $\supseteq$, let
$Qz \in Q\sX \cap Q\sY$ be given.  Then $Qz \in Q\sX$ so $z \in Q^{ -1
}Q\sX = \sX$ where we used $\sU \subseteq \sX$ and Remark
\ref{rmk:Verdier}.  Similarly we get $z \in \sY$ so $z \in \sX \cap \sY$ whence $Qz \in Q( \sX \cap \sY )$.

(i.b): The inclusion $\subseteq$ is clear.  To see $\supseteq$, let
$Qz \in Q\sX * Q\sY$ be given.  There is a distinguished triangle $Qx
\rightarrow Qz \rightarrow Qy$ in $\sT / \sU$ with $x \in \sX$, $y \in
\sY$.  By \cite[def.\ 8.4]{M} the triangle has the form $Q( a
\rightarrow c \rightarrow b )$ for a distinguished triangle $a
\rightarrow c \rightarrow b$ in $\sT$.  Then $Qa \cong Qx$ so $Qa \in
Q\sX$, and as in the proof of (i.a) we get $a \in \sX$.  Similarly we
get $b \in \sY$ so $c \in \sX * \sY$ whence $Qz \cong Qc \in Q( \sX *
\sY )$.

(ii): The inclusion $\supseteq$ follows from
$Q(\sU * \sX * \sU ) = Q(\sU )  * Q(\sX )* Q(\sU )$ and 
$Q\sU = 0$; indeed, $Q$ is an exact functor and 
$\sU$ is the kernel of $Q$ by \cite[cor.\
II.2.2.11(a)]{V}.

To see $\subseteq$, let $z \in Q^{ -1 }Q\sX$ be given.  This means $Qz
\in Q\sX$ so there is an isomorphism $Qx \stackrel{ \sim }{
  \rightarrow } Qz$ with $x \in \sX$.  By \cite[def.\ 2.1.11]{N}, the
isomorphism is represented by a diagram $x \stackrel{ \psi }{
  \leftarrow } w \stackrel{ \chi }{ \rightarrow } z$ in $\sT$ with
$\cone( \psi ) \in \sU$.  That is, the isomorphism equals $( Q\chi )(
Q\psi )^{ -1 }$.  Since $Q\psi$ and $( Q\chi )( Q\psi )^{ -1 }$ are
isomorphisms, so is $Q\chi$, and it follows that $0 \cong \cone( Q\chi
) \cong Q( \cone\, \chi )$ whence $\cone( \chi ) \in \sU$. From $x \in
\sX$ and $\cone( \psi ) \in \sU$ follows $w \in \sU * \sX$, and
together with $\cone( \chi ) \in \sU$ this shows $z \in \sU * \sX *
\sU$.
\end{proof}

\section{Proof of Theorem B}
\label{sec:proof_of_B}

Before the proof, recall again that a stable t-structure in $\sT$ is a
pair $( \sX,\sY )$ of subcategories which are stable under $\Sigma$
and satisfy $\Hom_{ \sT }( \sX,\sY ) = 0$ and $\sX * \sY = \sT$.

{\em Proof} of Theorem B, the {\em first} row of the table, in which
we have triangulated subcategories $\sX, \sY \subseteq \sT$ satisfying
$\sX * \sY = \sT$.  We set $\sU = \sX \cap \sY$ and consider the
quotient functor $Q : \sT \rightarrow \sT / \sU$.

We must show that $( Q\sX,Q\sY )$ is a stable t-structure in $\sT /
\sU$.

Since $\sX, \sY$ are closed under $\Sigma$ in $\sT$, it follows that
$Q\sX, Q\sY$ are closed under $\Sigma$ in $\sT/\sU$.  Moreover, $\sX *
\sY = \sT$ clearly implies $Q\sX * Q\sY = \sT / \sU$.  Finally,
Theorem A implies $\Hom_{ \sT / \sU }( Q\sX , Q\sY ) = 0$.
$\;\;\;\;$ \hfill $\Box$
\medskip

{\em Proof} of Theorem B, the {\em second} row of the table, in which
we have triangulated subcategories $\sX,\sY,\sZ \subseteq \sT$
satisfying $\sX * \sY = \sY * \sZ = \sT$.  We set $\sU = \langle \sX
\cap \sY , \sY \cap \sZ \rangle $ and consider the quotient functor $Q
: \sT \rightarrow \sT / \sU$.

We must show that the pairs $( Q\sX,Q\sY )$ and $( Q\sY,Q\sZ )$ are
stable t-structures in $\sT / \sU$.  It suffices to show a proof for
the first pair as the second is similar.

Set
\[
  \sY' = \sY \cap \sZ
  \;\;,\;\; 
  \widetilde{ \sX } = \langle \sX,\sY' \rangle.
\]
Since $\sY' \subseteq \sU$ we have $Q\widetilde{ \sX } =
Q\sX$, so it is enough to show that $( Q\widetilde{ \sX } , Q\sY )$ is
a stable t-structure.

Since $\widetilde{ \sX }, \sY$ are closed under $\Sigma$ in $\sT$, it
follows that $Q\widetilde{ \sX }, Q\sY$ are closed under $\Sigma$ in
$\sT/\sU$.  Moreover, $\sX * \sY = \sT$ implies $\widetilde{ \sX } *
\sY = \sT$ whence $Q\widetilde{ \sX } * Q\sY = \sT / \sU$ is clear.
Finally, we have $\sY' \subseteq \sY$ whence Lemma
\ref{lem:angles}(ii) gives the second ``$=$'' in
\[
  \widetilde{ \sX } \cap \sY
  = \langle \sX,\sY' \rangle \cap \sY
  = \langle \sX \cap \sY , \sY' \rangle
  = \langle \sX \cap \sY , \sY \cap \sZ \rangle
  = \sU.
\]
Since $\widetilde{ \sX } * \sY = \sT$ is triangulated, Theorem A
implies $\Hom_{ \sT / \sU }( Q\widetilde{ \sX },Q\sY ) = 0$. 
\hfill $\Box$
\medskip

{\em Proof} of Theorem B, the {\em third} row of the table, in which
we have triangulated subcategories $\sX,\sY,\sZ \subseteq \sT$
satisfying $\sX * \sY = \sY * \sZ = \sZ * \sX = \sT$.  We set $\sU =
\langle \sX \cap \sY , \sY \cap \sZ , \sZ \cap \sX \rangle $ and
consider the quotient functor $Q : \sT \rightarrow \sT / \sU$.

We must show that the pairs $( Q\sX,Q\sY )$, $( Q\sY,Q\sZ )$, $(
Q\sZ,Q\sX )$ are stable t-structures in $\sT / \sU$.  It suffices to
show a proof for the first pair as the others are similar.

Indeed, this can be shown analogously to the previous case by setting 
\begin{align*}
  & \sY' = \sY \cap \sZ
  \;\;,\;\;
  \widetilde{ \sX } = \langle \sX,\sY' \rangle \\
  & \sX' = \sZ \cap \sX
  \;\;,\;\;
  \widetilde{ \sY } = \langle \sY,\sX' \rangle.
\end{align*}
Then $Q\widetilde{ \sX } = Q\sX, Q\widetilde{ \sY } = Q\sY$, and it is
enough to show that $( Q\widetilde{ \sX } , Q\widetilde{ \sY } )$ is a
stable t-structure.

As above, $Q\widetilde{ \sX }, Q\widetilde{ \sY }$ are closed under
$\Sigma$ in $\sT/\sU$ and $Q\widetilde{ \sX } * Q\widetilde{ \sY } =
\sT / \sU$.  Finally, since $\sX' \subseteq \sX, \sY' \subseteq \sY$,
a computation based on Lemma \ref{lem:angles}(i and ii) gives
$\widetilde{ \sX } \cap \widetilde{ \sY } = \sU$.  As above, Theorem A
then implies $\Hom_{ \sT / \sU }( Q\widetilde{ \sX },Q\widetilde{ \sY
} ) = 0$.
\hfill $\Box$

\section{Proof of Theorem C}
\label{sec:proof_of_C}

When $\sX,\sY \subseteq \sT$ are thick subcategories satisfying $\sX
\cap \sY = \sU$ and $\sX * \sY = \sT$, the first item of the table in
Theorem B shows that $( Q\sX , Q\sY )$ is a stable t-structure in $\sT
/ \sU$.  Hence $( \sX,\sY ) \mapsto ( Q\sX,Q\sY )$ is indeed a map
\begin{equation}
\label{equ:map}
  \bigg\{\, ( \sX , \sY )
  \,\bigg|
    \begin{array}{l}
      \mbox{ $\sX, \sY \subseteq \sT$ thick subcategories } \\[2mm]
      \mbox{ with $\sX \cap \sY = \sU$ and $\sX * \sY = \sT$ }
    \end{array}
  \,\bigg\}
  \rightarrow
  \{\, \mbox{ stable t-structures in $\sT/\sU$ } \,\},
\end{equation}
as claimed in Theorem C.

The map \eqref{equ:map} is injective because of Remark
\ref{rmk:Verdier}.

The map \eqref{equ:map} is surjective: A stable t-structure in $\sT /
\sU$ consists of two thick subcategories.  By Remark \ref{rmk:Verdier}
it has the form $( Q\sX,Q\sY )$ for two unique thick subcategories
$\sX,\sY$ satisfying $\sU \subseteq \sX \cap \sY \subseteq \sT$.  To
complete the proof, we must show $\sX \cap \sY = \sU$ and $\sX * \sY =
\sT$.

Lemma \ref{lem:Q}(i.a) gives $Q( \sX \cap \sY ) = Q\sX \cap Q\sY =
0$.  This implies $\sX \cap \sY = \sU$ by Remark \ref{rmk:Verdier}.

Note that as a consequence,
\[
  \Hom_{ \sT / \sX \cap \sY }( Q\sX,Q\sY )
  = \Hom_{ \sT / \sU }( Q\sX,Q\sY ) = 0
\]
so $\sX * \sY$ is triangulated by Theorem A.

Lemma \ref{lem:Q}(i.b) gives $Q( \sX * \sY ) = Q\sX * Q\sY = \sT / \sU
= Q\sT$.  This implies $\sX * \sY = \sT$ by Remark
\ref{rmk:Verdier} since $\sX * \sY$ is a triangulated subcategory
satisfying $\sU \subseteq \sX * \sY \subseteq \sT$.  \hfill $\Box$

\section{Another result on subcategories of extensions}
\label{sec:another}

\begin{Theorem}
Let $\sX, \sY, \sV$ be thick subcategories of $\sT$ satisfying $\sX
\cap \sY \subseteq \sV \subseteq \sT$ and $\sX * \sY = \sT$.  The
following are equivalent.
\begin{enumerate}

  \item  $\sX * \sV$ is triangulated.

\smallskip

  \item  $\sV * \sY$ is triangulated.

\smallskip

  \item  $\sV = ( \sX \cap \sV ) * ( \sY \cap \sV )$. 

\end{enumerate}
\end{Theorem}

\begin{proof}
Write $\sU = \sX \cap \sY$ and let $Q : \sT \rightarrow \sT / \sU$ be
the quotient functor.  Note that $Q\sX$, $Q\sY$, $Q\sV$ are thick (and
in particular triangulated) subcategories of $\sT / \sU$ by Remark
\ref{rmk:Verdier}.

We claim that each condition in the theorem is equivalent to the
corresponding condition in the following list.
\renewcommand{\labelenumi}{(\roman{enumi}')}
\begin{enumerate}

  \item  $Q\sX * Q\sV$ is triangulated.

\smallskip

  \item  $Q\sV * Q\sY$ is triangulated.

\smallskip

  \item  $Q\sV = ( Q\sX \cap Q\sV ) * ( Q\sY \cap Q\sV )$.

\end{enumerate}
\renewcommand{\labelenumi}{(\roman{enumi})}

(i)$\Leftrightarrow$(i'): By Lemma \ref{lem:factorization} this amounts to
\[
  \sV * \sX \subseteq \sX * \sV
  \Leftrightarrow
  Q\sV * Q\sX \subseteq Q\sX * Q\sV.
\]
But $\Rightarrow$ follows from Lemma \ref{lem:Q}(i.b) and
$\Leftarrow$ can be seen as follows using Lemma \ref{lem:Q}(i.b and
ii).
\begin{align*}
  Q\sV * Q\sX \subseteq Q\sX * Q\sV
  & \Leftrightarrow Q( \sV * \sX ) \subseteq Q( \sX * \sV ) \\ 
  & \Rightarrow Q^{-1}Q( \sV * \sX ) \subseteq Q^{-1}Q( \sX * \sV ) \\
  & \Rightarrow \sU * \sV * \sX * \sU \subseteq \sU * \sX * \sV * \sU \\
  & \Leftrightarrow \sV * \sX \subseteq \sX * \sV.
\end{align*}
The last $\Leftrightarrow$ holds since $\sU \subseteq \sV$ and $\sU
\subseteq \sX$. 

(ii)$\Leftrightarrow$(ii') and (iii)$\Leftrightarrow$(iii'): Similar.


We now show that (i') through (iii') are equivalent.

(i')$\Rightarrow$(iii'): The inclusion $\supseteq$ is obvious. We shall show $\subseteq$. 
The condition $\sX * \sY = \sT$ clearly
implies $Q\sX * Q\sY = \sT / \sU$ so given $Qv \in Q\sV$
there is a distinguished triangle
\begin{equation}
\label{equ:triangle}
  Qx \rightarrow Qv \rightarrow Qy \rightarrow \Sigma Qx
\end{equation}
in $\sT / \sU$ with $Qx \in Q\sX$, $Qy \in Q\sY$.  This shows $Qy \in
Q\sV * Q\sX$ but (i') means $Q\sV * Q\sX \subseteq Q\sX * Q\sV$ by
Lemma \ref{lem:factorization}, so $Qy \in Q\sX * Q\sV$.  Hence there is a distinguished
triangle $Qx_1 \rightarrow Qy \rightarrow Qv_1$ with $Qx_1 \in Q\sX$,
$Qv_1 \in Q\sV$, but $\Hom_{ \sT / \sU }( Q\sX,Q\sY ) = 0$ by Theorem
A so the triangle splits whence $Qy$ is a direct summand of $Qv_1$.
Since $Q\sV$ is thick this implies $Qy \in Q\sV$.  The distinguished
triangle \eqref{equ:triangle} then shows $Qx \in Q\sV$.  We have shown
$Qx \in Q\sX \cap Q\sV$, $Qy \in Q\sY \cap Q\sV$ so
\eqref{equ:triangle} shows (iii').

(iii')$\Rightarrow$(i'): When (iii') holds we have
\[
  Q\sX * Q\sV
  = Q\sX * \big( ( Q\sX \cap Q\sV ) * ( Q\sY \cap Q\sV ) \big)
  = \big( Q\sX * ( Q\sX \cap Q\sV ) \big) * ( Q\sY \cap Q\sV )
  = Q\sX * ( Q\sY \cap Q\sV )
\]
so it is enough to see that $Q\sX * ( Q\sY \cap Q\sV )$ is
triangulated.  However, $Q\sX$ and $Q\sY \cap Q\sV$ are both
triangulated by Remark \ref{rmk:Verdier} and $\Hom_{ \sT/\sU }(
Q\sX,Q\sY \cap Q\sV ) = 0$ because $\Hom_{ \sT/\sU }( Q\sX,Q\sY ) = 0$
by Theorem A.  Hence $Q\sX * ( Q\sY \cap Q\sV )$ is triangulated by
Corollary \ref{cor:A}.

(ii')$\Leftrightarrow$(iii') follows by similar arguments.
\end{proof}

\section{Examples}
\label{sec:examples}

\subsection{The homotopy category of projective modules}

Let $R$ be an Iwa\-na\-ga-Go\-ren\-stein ring, that is, a
noetherian ring which has finite injective dimension from either side
as a module over itself.  Let $\sT = \sK_{ ( \b ) }( \Prj\, R )$ be
the homotopy category of complexes of projective right-$R$-modules with
bounded homology.  Define subcategories of $\sT$ by
\[
  \sX = \sK^{ - }_{ (\b) }( \Prj\, R )
  \;\;,\;\;
  \sY = \sK_{ \ac }( \Prj\, R )
  \;\;,\;\;
  \sZ = \sK^{ + }_{ (\b) }( \Prj\, R )
\]
where $\sK^{ - }_{ (\b) }( \Prj\, R )$ is the isomorphism closure of
the class of complexes $P$ with $P^i = 0$ for $i \gg 0$ and $\sK^{ +
}_{ (\b) }( \Prj\, R )$ is defined analogously with $i \ll 0$, while
$\sK_{ \ac }( \Prj\, R )$ is the subcategory of acyclic (that is,
exact) complexes.

Note that $\sY$ is equal to $\sK_{ \tac }( \Prj\, R )$, the
subcategory of totally acyclic complexes, that is, acyclic complexes
which stay acyclic under the functor $\Hom_R( - , Q )$ when $Q$ is
projective, see \cite[cor.\ 5.5 and par.\ 5.12]{IK}.

By \cite[prop.\ 2.3(1), lem.\ 5.6(1), and rmk.\ 5.14]{IKM} there are
stable t-structures $( \sX , \sY )$, $( \sY , \sZ )$ in $\sT$, so
\[
  \sX * \sY = \sY * \sZ = \sT
  \;\;,\;\;
  \sX \cap \sY = \sY \cap \sZ = 0.
\]

If $P \in \sT$ is given, then there is a distinguished triangle $P^{
\geq 0 } \rightarrow P \rightarrow P^{ < 0 }$ where $P^{ \geq 0 }$ and
$P^{ < 0 }$ are the hard truncations.  Since $P^{ \geq 0 } \in \sZ$ and
$P^{ < 0 } \in \sX$, we have
\[
  \sZ * \sX = \sT.
\]

Finally,
\[
  \sZ \cap \sX
  =
  \sK^{ + }_{ (\b) }( \Prj\, R ) \cap \sK^{ - }_{ (\b) }( \Prj\, R )
  =
  \sK^{ \b }( \Prj\, R )
\]
is the isomorphism closure of the class of bounded complexes.  If we
use an obvious shorthand for quotient categories, then the last row of
the table in Theorem B now provides a triangle of recollements
\[
  \big( 
  \sK^{ - }_{ (\b) } / \sK^{ \b }( \Prj\, R )
  \;,\;
  \sK_{ \ac }( \Prj\, R )
  \;,\;
  \sK^{ + }_{ (\b) } / \sK^{ \b }( \Prj\, R )
  \big)
\]
in $\sK_{ (\b) } / \sK^{ \b }( \Prj\, R )$.  The reason that we can
write $\sK_{ \ac }( \Prj\, R )$ instead of its image in $\sK_{
  (\b) } / \sK^{ \b }( \Prj\, R )$ is that the two are equivalent by
\cite[prop.\ 1.5]{IKM}.

This example and its analogue with finitely generated projective
modules were first obtained in \cite[thms.\ 2.8 and 5.8]{IKM} and
motivated the definition of triangles of recollements.

\subsection{The symmetric Auslander category}

Let $R$ be a noetherian commutative ring with a dualizing complex.  In
\cite[def.\ 2.1]{JK} we introduced the symmetric Auslander
ca\-te\-go\-ry $\sA^{\s}(R)$ as the isomorphism closure in $\sK_{ (\b)
}( \Prj\, R )$ of the class of complexes $P$ for
which there exist $T, U \in \sK_{ \tac }( \Prj\, R )$ such that $P^{
  \ll 0} = T^{\ll 0}$ and $P^{\gg 0} = U^{\gg 0}$.

The methods of \cite[sec.\ 2]{JK} show that there are stable
t-structures
\[
  \big( \sA(R) , \sK_{ \tac }( \Prj\, R ) \big)
  \;\;,\;\;
  \big( \sK_{ \tac }( \Prj\, R ) , S( \sB(R) ) \big)
\]
in $\sA^{ \s }( R )$ with $S( \sB(R) ) * \sA( R ) = \sA^{ \s }( R )$.
Here $\sA(R)$ and $\sB(R)$ are the Auslander and Bass categories of
$R$, see \cite[sec.\ 3]{AF}, and $S$ is the functor introduced in
\cite[sec.\ 4.3]{IK}.  
Since $\sA(R)$ (resp. $S( \sB (R))$) consists of right-bounded (resp. left-bounded) complexes 
$P$ 
such that $P^{
  \ll 0} = T^{\ll 0}$ (resp. $P^{\gg 0} = T^{\gg 0}$ )
for some $T  \in \sK_{ \tac }( \Prj\, R )$  \cite[Rem. 2.2]{JK}, 
we have $\sA( R ) \cap S( \sB(
R ) ) = \sK^{ \b }( \Prj\, R )$. 

With $\sT = \sA^{ \s }( R )$ and
\[
  \sX = \sA(R)
  \;\;,\;\;
  \sY = \sK_{ \tac }( \Prj\, R )
  \;\;,\;\;
  \sZ = S( \sB(R) ),
\]
the last row of the table in Theorem B recovers the triangle of
recollements 
\[
  \big(
  \sA(R) / \sK^{ \b }( \Prj\, R )
  \;,\;
  \sK_{ \tac }( \Prj\, R )
  \;,\;
  S( \sB(R) ) / \sK^{ \b }( \Prj\, R )
  \big)
\]
in $\sA^{ \s }( R ) / \sK^{ \b }( \Prj\, R )$ first obtained in
\cite[thm.\ 2.10]{JK}.  As in the previous example, the category
$\sK_{ \tac }( \Prj\, R )$ is equivalent to its image in $\sA^{
  \s }( R ) / \sK^{ \b }( \Prj\, R )$ so we can write $\sK_{ \tac }(
\Prj\, R )$ instead of the image.

\bigskip
{\bf Acknowledgement. }
We gratefully acknowledge support from JSPS Grant 22540053 and a
Scheme 2 Grant from the London Mathematical Society.

\end{document}